\documentclass[preprint,12pt]{elsarticle}

\usepackage{amssymb}
\usepackage{amsmath}
\usepackage{bm}
\DeclareMathOperator*{\argmin}{arg\,min}
\renewcommand{\vec}{\bm}

\usepackage[linesnumbered,lined,boxed,commentsnumbered]{algorithm2e}

\usepackage{enumitem}
\usepackage{subcaption}

\usepackage{hyperref}

\usepackage{nomencl}
\makenomenclature
\renewcommand\nomgroup[1]{%
  \item[\emph{
  \ifstrequal{#1}{A}{Acronyms}{%
  \ifstrequal{#1}{B}{Roman symbols}{%
  \ifstrequal{#1}{C}{Greek symbols}{%
  \ifstrequal{#1}{E}{Superscripts}{
  \ifstrequal{#1}{D}{Subscripts}{
  \ifstrequal{#1}{O}{Other symbols}{%
  }}}}}}%
}]}

\newcommand{\nomRoman}[1][]{\nomenclature[B,#1]}
\newcommand{\nomGreek}[1][]{\nomenclature[C,#1]}
\newcommand{\nomSuper}[1][]{\nomenclature[E,#1]}
\newcommand{\nomSub}[1][]{\nomenclature[D,#1]}
\newcommand{\nomAcro}[1][]{\nomenclature[A,#1]}

\journal{Computers and Chemical Engineering}

\begin{document}

\begin{frontmatter}



\title{Piecewise linear approximation for MILP leveraging piecewise convexity to improve performance}


\author[IET]{Felix Birkelbach}
\author[IET]{David Huber}
\author[IET]{René Hofmann}

\affiliation[IET]{organization={Institute for Energy Systems and Thermodynamics},
            addressline={Getreidemarkt 6/E302}, 
            city={Vienna},
            postcode={1060}, 
            state={},
            country={Austria}}



\begin{abstract}
To realize adaptive operation planning with MILP unit commitment, piecewise-linear approximations of the functions that describe the operating behavior of devices in the energy system have to be computed. We present an algorithm to compute a piecewise-linear approximation of a multi-variate non-linear function. The algorithm splits the domain into two regions and approximates each region with a set of hyperplanes that can be translated to a convex set of constraints in MILP.  
The main advantage of this ``piecewise-convex approximation'' (PwCA) compared to more general piecewise-linear approximation with simplices is that the MILP representation of PwCA requires only one auxiliary binary variable.  For this reason, PwCA yields significantly faster solving times in large MILP problems where the MILP representation of certain functions has to be replicated many times, such as in unit commitment.
To quantify the impact on solving time, we compare the performance using PwCA with the performance of simplex approximation with logarithmic formulation and show that PCA outperforms the latter by a big margin. For this reason, we conclude that PCA will be a useful tool to set up and solve large MILP problems such as arise in unit commitment and similar engineering optimization problems.

\end{abstract}

\begin{graphicalabstract}
\end{graphicalabstract}

\begin{highlights}
\item We present an algorithm to compute a piecewise-linear approximation of a data set
\item We discuss the accuracy of the approximation and show that it yields fast solving times in MILP
\item We show that piecewise-convex approximation is a useful tool for engineering optimization problems 
\end{highlights}

\begin{keyword}
MILP \sep piecewise-linear approximation \sep unit commitment
\end{keyword}

\end{frontmatter}



\section{Introduction}
\label{sec:intro}

%
%
%


MILP has been applied successfully in a wide range of engineering applications such as combined design- and operation optimization \cite{Halmschlager2022}, process control \cite{Fuhrmann2022} and unit commitment \cite{halmschlager_optimizing_2021}, to name just a few. MILP solvers have made big advances in recent years, which allow to solve ever larger problems in less time \cite{Koch2022}. The availability of mature solvers makes MILP interesting for industrial applications. 

Considering that optimization problems in engineering are generally non-linear, mixed integer non-linear programming (MINLP) approaches would be the natural choice to solve these problems. Even though there have been big advances in the area of deterministic MINLP solvers recently, they are still limited to comparably small problem instances \cite{Kronqvist2019}. Meta-heuristic solvers for MINLP are not an option for many applications, because they cannot guarantee to find the global optimal solution. Another option is to approximate the non-linear parts of the problem and use mixed integer linear programming (MILP). Even though there will be some approximation error, this drawback is generally out-weighted by the availabiltity of mature solvers for MILP problems. The central challenge is to find viable (piecewise) linear approximations of the non-linear functions \cite{brito2020}.


Our goal is to develop a framework for continuous operation planning for energy systems \cite{schwarzmayr_development_2022}, where the models that describe the operating behavior of devices in the energy system are regularly updated based on historic data to account for gradual changes of the operating behavior. The operation planning will be done with a MILP unit commitment problem. The operating behavior of each device is represented by a non-linear function. To realize the adaptive modeling of the operating behavior, we need to compute piecewise-linear approximations of the data set that represents the non-linear operating behavior. 

Finding an efficient piecewise-linear approximation is especially crucial for unit commitment and similar optimization problems, where the MILP representation of certain non-linear functions has to be replicated many times. In unit commitment problems, the non-linear operating behavior of certain devices has to be replicated at each time step in the MILP model. A more efficient MILP representation of the operating behavior can have a significant impact on the solving time of the MILP problem.

We are interested in piecewise-linear approximation methods that can be applied to continuous functions with two or more independent variables. The approximation must be able to capture key features of the non-linear function such as directions of positive and negative curvature and it must have an efficient (in terms of impact on the solving time) MILP representation. 

An interesting approach was used by Koller et~al.~\cite{Koller_2019}: they divided a non-linear function that described the operating behavior of a thermal storage into two regions and set up a convex approximation of the function in each region so that the approximation was continuous at the interface. This resulted in an excellent approximation of their function and only required one auxiliary binary variable in the MILP problem to toggle between the two regions. Even though, this approach puts certain limitations on the shape of the non-linear function (i.e. two regions where the function is approximately convex), it is considerably more versatile than simple convex approximation and, at the same time, should yield very fast solving times in MILP. Koller et~al.\ set up their approximation by hand, which is why their approach cannot be readily transferred to other applications. Though, if it was possible to set up this type of approximation based on a data set that describes the non-linear function, these models could be very useful for MILP applications where performance is critical and certain non-linear features cannot be neglected, such as our adaptive operation planning framework.


In this paper we present a novel algorithm to approximate a multi-variate continuous non-linear function by splitting it in two regions and computing a convex approximation for each region. In the remainder of this paper, we will refer to this approach as piecewise-convex approximation (PwCA). In the next section we discuss other methods for piecewise linear approximation methods and their MILP representations. In Section~\ref{sec:algorithm} we describe the approximation algorithm and how a PwCA can be translated to MILP. In Section~\ref{sec:results} we assess the performance of PwCA. We show that PwCA is a considerably better approximation than simple convex approximation for many non-linear functions and that it outperforms simplex approximation in large MILP optimization problems. From this, we conclude in Section~\ref{sec:conclusion} that PwCA fills an important gap among piecewise linear approximation methods for MILP and that it is a useful tool for solving large unit commitment problems and similar engineering optimization problems where the MILP representation of a non-linear function has to be replicated many times.

\section{Background}

The simplest, albeit also quite limited, method to linearize a non-linear function is to compute a linear approximation. If this approximation is sufficiently accurate, it is by a big margin the most efficient way to linearize functions for MILP. If it is not sufficiently accurate, two more sophisticated methods can be applied: piecewise linear approximation with simplices and convex/concave approximation with (hyper-)planes.

Piecewise linear approximation with simplices is the most universal method to incorporate a non-linear function in MILP. The domain of the function is triangulated and the function is interpolated on each simplex of the triangulation. (For functions with one independent variable, this is equivalent to linear spline approximation.)  Any surjective function, even highly non-linear ones, can be approximated with arbitrary accuracy (i.e. approximation error).  By increasing the number of linear segments, the accuracy of the approximation can be increased. Though, this accuracy comes at a cost: the more segments a function is divided into, the more constraints and auxiliary variables are introduced into the optimization problem, which can quickly make it intractable. 
For functions with one independent variable, the simplex approximation can be translated to MILP as a special ordered set of type 2 (SOS2). This is highly efficient and covered in introductory text books on MILP. Functions with more than one independent variables, cannot be translated as efficiently \cite{Vielma_2010}. Especially for multi-variate functions the number of auxiliary variables increases quickly. Various MILP formulations are available, which differ in the number of auxiliary variables and constraints \cite{Vielma_2010}. Though, even with specialized MILP formulations optimized for a large number of simplices \cite{Vielma_2011}, the maximum accuracy which is still tractable may be severely limited.
%

With convex/concave approximation, the non-linear function is incorporated into the MILP problem by super-positioning multiple linear constraints. Thus, the non-linear function has to be approximated by a set of hyperplanes. As the name suggests, this approximation will only yield good results, if the non-linear function is either convex or concave. For many engineering problems, where the non-linear function is neither convex or concave, the approximation will be poor.
The main advantage is that the MILP representation is very efficient. For representing convex approximations in minimization problems or concave approximations for maximization problems, not a single auxiliary variables is required. For concave approximations in minimization problems and convex approximations in maximization problems, one binary variable per constraint has to be introduced. Overall, this results in considerably less auxiliary variables than with simplex approximation, and consequently these approximations generally yield very good MILP performance. Another advantage is that this approach extends to functions with many independent variables with very little overhead.  
For simplicity, we will refer to this approach as convex approximation in the remainder of this paper, even though it can also be applied to concave functions.

Many non-linear functions in engineering problems do not have an exceedingly complicated shape. If convex approximation is sufficiently accurate, it is the natural choice because it is easy to set up and performs very well. Though, convex approximation will fail to capture key features of the function, if the function has both a direction with positive and a direction with negative curvature. Then, the only other option is piecewise linear approximation with simplices, which introduces a large number of auxiliary binary variables and consequently has a negative impact on performance ---especially for multi-variate functions.

We propose piecewise-convex approximation (PwCA), which can capture features of non-linear function that simple convex approximation cannot. If the accuracy of the PwCA is sufficient, it yields much better MILP performance than simplex approximation because the MILP representation of PwCA requires only one auxiliary binary variable. 

Finding a good PwCA for a given non-linear function is surprisingly challenging ---\,and doing it with an automated algorithm even more so. Not only do we need to find the right position to split the function into two regions so that a convex approximation of the function in each region is accurate (in terms of deviation from the non-linear function), but also do we need to set up the convex approximations in both regions so that the hyperplanes coincide on the interface. Otherwise the approximation will not be continuous, which may cause problems when translating the approximation to MILP. Na\"ively extending the fitting approach of simple convex approximation to approximation with two regions and implementing the continuity requirement as a constraint in the fitting procedure (i.e. treating the fitting problem as a constrained non-linear optimization problem) does not work. The continuity constraint is so complex that it inhibits convergence. A different approach is required to compute a PwCA.

\section{The algorithm}
\label{sec:algorithm}

Consider a continuous function $y = f(x_1,\ldots,x_{n-1})$ with $n-1$ independent variables in the domain $x_{j,\mathrm{min}} \le x_j \le x_{j,\mathrm{max}}$. This function describes a hyper-surface in $n$-dimensional space. We assume that the non-linear function is given by a set of $m\in \{1,\ldots, N_\mathrm{data}\}$  data points $\{x'_{1,m},\ldots,x'_{n-1,m}, y'_m \}$. This data set may be derived from evaluating the non-linear function on a grid, simulation data or even experimental data. The goal of the algorithm is to compute an accurate piecewise linear approximation of this non-linear function, i.e. to minimize the deviation between the estimate $\hat{y}$ and the data $y_m$ by adjusting the model parameters. 

Piecewise-convex approximation (PwCA) can be seen as an extension of simple convex approximation. For this reason, we will first describe how to compute a simple convex approximation and then we will move to PwCA. 

A Matlab implementation of these algorithms is available at \url{https://gitlab.tuwien.ac.at/iet/public/milptools}.

\subsection{Computing a simple convex approximation}
\label{sec:fitSC}

 With convex approximation, the function is approximated by a set of linear functions, i.e. hyperplanes, that can be translated to a convex set of constraints in MILP. Assume that each hyperplane $i \in \{1,\ldots,N_\mathrm{hyp} \}$ is given by a set of $n+1$ coefficients $\vec{a}_i$ so that $[1,x_1,\ldots,x_{n-1},y] \vec{a}_i = 0$ for all points on the hyperplane. Then the estimate $\hat{y}(\vec{x})$ is given by
\begin{equation}
    \hat{y}(\vec{x}) = \max \{y_i(\vec{x}) : i \in \{1,\ldots,N_\mathrm{hyp}\} \} 
    \label{eq:yhat}
\end{equation}
where
\begin{equation}
     y_i(\vec{x}) = \left( a_{i,0} + \sum_{j=1}^{n-1} a_{i,j} x_j \right) \frac{1}{a_{i,n}}
     \label{eq:yi}
\end{equation}
is the $y$-value of the $i$th hyperplane evaluated at $\vec{x} = [x_1,\ldots,x_{n-1}]$.

The objective of the fitting algorithm is to minimize the deviation in terms of the sum of squared errors
\begin{equation}
    \argmin_{\vec{a}_1,\ldots,\vec{a}_{N_\mathrm{hyp}}} \sum_{m=1}^{N_\mathrm{data}}\big(\hat{y}(x'_{1,m},\ldots,x'_{n-1,m}) - y'_m \big)^2
\end{equation}
by adjusting the hyperplane parameters $\vec{a}_1$ through $\vec{a}_{N_\mathrm{hyp}}$. 

This is an unconstrained non-linear minimization problem. It turns out to converge remarkably well, regardless of the initial guess for the hyperplane parameters $\vec{a}_i$. Only if many hyperplanes are used for the approximation, it can happen that some hyperplanes are ``unused'', i.e. that they are dominated by the other hyperplanes on the whole domain. This can be remedied by penalizing the distance between $y_i$ at the corners of the domain and $\max \{y'_m\}$, which drives hyperplanes upwards until they contribute to the approximation.

The same approach can be used to fit concave functions. In Eq.~\eqref{eq:yhat}, $\max y_i$ has to be replaced by $\min y_i$ or the identity $\min y_i = -\max{-y_i}$ can be used to transform the problem.

\subsection{Computing a piecewise-convex approximation}
\label{sec:fitDC}

With PwCA, the non-linear function is divided into two regions and for each region, a convex approximation is computed. The two convex approximations have to coincide at the interface that separates the two regions. Otherwise the PwCA would not be continuous, which may result in unwanted effects in the MILP problem. To ensure continuity, PwCA always uses an even number of hyperplanes $N_\mathrm{hyp}$. Each pair of planes exactly intersects at the interface.

The na\"ive approach of extending the algorithm for simple convex approximation to PwCA would be to implement the continuity requirement at the interface as a constraint. Unfortunately, it turns out that the resulting constrained non-linear optimization problem does not converge because the continuity constraint is so complex that it inhibits convergence. The problem has to be reformulated so that the continuity is inherent in the way that the approximation is set-up. The goal is to parameterize the approximation model in a way that we again arrive at an unconstrained problem.

We propose to use a series of rotations and shifts to set up the approximation model. Instead of the hyperplane parameters $\vec{a}_i$ we use a set of rotation angles and distances to define the position of the hyperplanes. The clever thing with this approach is that the hyperplanes are defined via the intersection at the interface. This guarantees continuity at the interfaces and we again end up with an unconstrained non-linear optimization problem, which converges reliably. 

\begin{figure}
    \centering
    \makebox[\textwidth][c]{
            \includegraphics{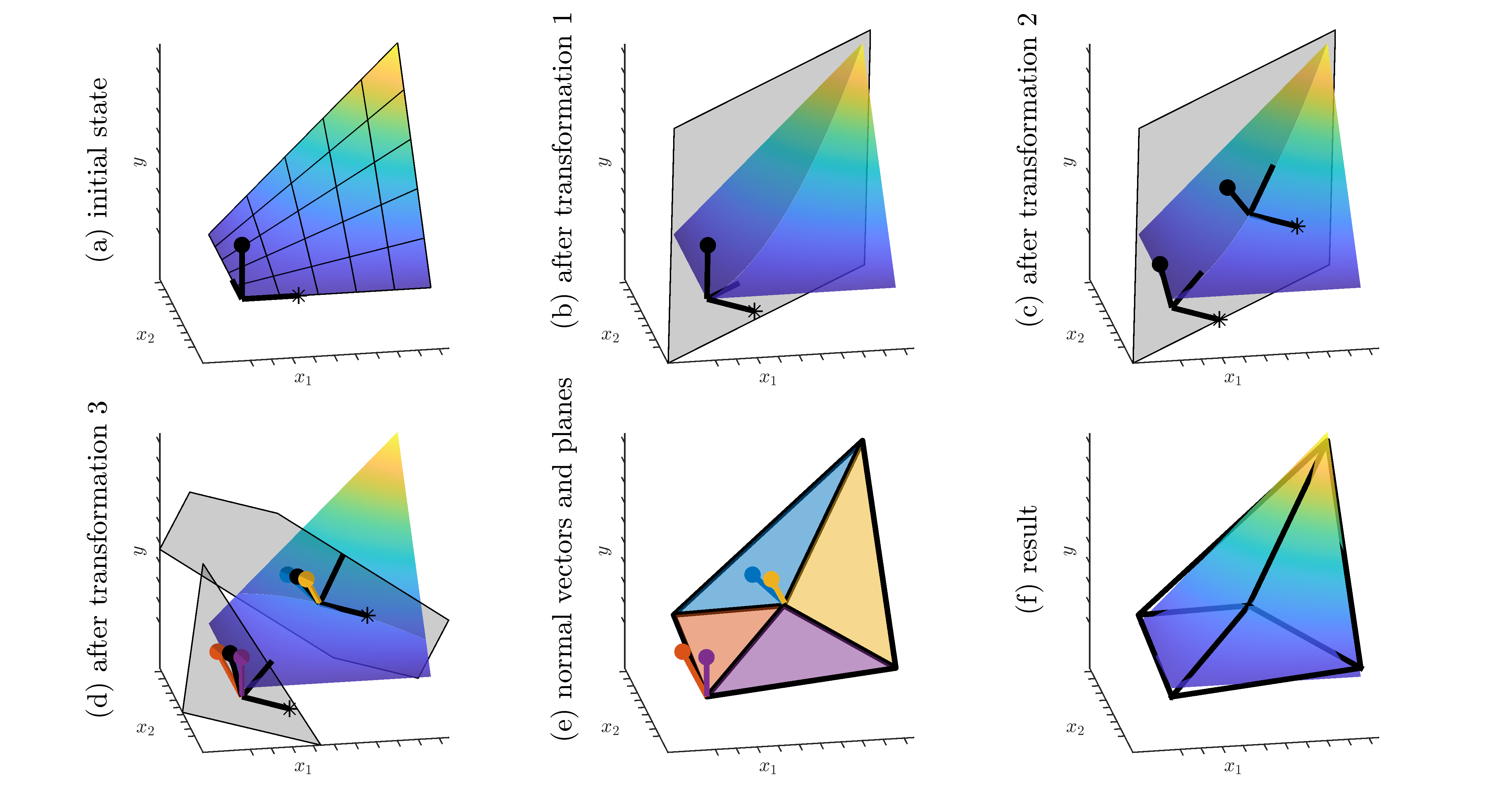}
    }
        \caption{Illustration of the steps to set up a PwCA with two pairs of hyperplanes. The vector with the dot is the $\vec{y}$ vector, which is the normal vector of the hyperplanes. The vector with the asterisks is the normal vector of the interface.}
        \label{fig:doubleConvex}
\end{figure}

For the 3-dimensional case, the rotations and shifts to set up the model are illustrated in Figure~\ref{fig:doubleConvex}. In 3D, hyperplanes are ordinary planes and hyperlines are ordinary lines. Figure~\ref{fig:doubleConvex} shows a PwCA with two pairs of hyperplanes.

The starting point for the algorithm is the orthonormal basis with $n$ basis vectors in the directions of $\vec{x}_1$ to $\vec{x}_{n-1}$ and $\vec{y}$ respectively (see Figure~\ref{fig:doubleConvex}\,a). Each hyperplane, including the interface, is uniquely identified by a normal vector and a point on the hyperplane. After performing the rotations and shifts, the transformed $\vec{y}$ vector will be the normal vector of each hyperplane. In Figure~\ref{fig:doubleConvex} it is marked by a dot at the end. Another vector has to be selected as the normal vector of the interface. In principle, this can be any of the other base vectors. Here, we chose the $\vec{x}_1$ vector. In Figure~\ref{fig:doubleConvex} it is marked by an asterisks.

For the rotations, we consider basic rotations, i.e. rotations in each of the planes that is spanned by the pairwise combinations of the basis vectors. In two dimensions there is one basic rotation, in three dimensions three basic rotations and in the general $n$-dimensional case there are $\binom{n}{2}$ basic rotations.


Transformation~1 is the rotation of the initial basis and shifting it into the direction of $\vec{x}_1$. Then, the vectors $\vec{x}_2$ to $\vec{x}_{n-1}$ and $\vec{y}$ span the hyperplane that is the interface and $\vec{x}_1$ is the normal vector of the interface (marked by an asterisk in Figure~\ref{fig:doubleConvex}\,b).  Transformation~1 is characterized by $\binom{n}{2}$ rotation angles $\vec{r}_1$ and one shift distance $s_1$. It is the same for all hyperplanes of the model.

Transformation~2 is to rotate the basis again, but only in the planes spanned by the pairwise combination of the vectors $\vec{x}_2$ to $\vec{x}_{n-1}$ and $\vec{y}$. The orientation of the vector $\vec{x}_1$, the normal vector of the interface, does not change. Then, the basis is shifted into the direction of $\vec{y}$ to get the intersection hyperline of each pair of hyperplanes on the interface. The intersection hyperline is spanned by the vectors $\vec{x}_2$ to $\vec{x}_{n-1}$ and the vectors $\vec{x}_1$ and $\vec{y}$ are the normal vectors. 
In the 3D case, which is illustrated in Figure~\ref{fig:doubleConvex}, each intersection hyperline is a simple line, so only one vector $\vec{x}_2$ is required to define the intersection. Since the model in Figure~\ref{fig:doubleConvex} has two pairs of hyperplanes, two bases shown. 
To describe transformation~2, $\binom{n-1}{2}$ rotation angles $\vec{r}_{2,i}$ plus one shift distance $s_{2,i}$ is needed for each pair of hyperplanes $i \in \{1,\ldots,N_\mathrm{hyp}/2\} $.

Transformation~3 is to rotate the basis of the intersection in the plane that is spanned by $\vec{x}_1$ and $\vec{y}$. This rotation does not affect the vectors $\vec{x}_2$ to $\vec{x}_{n-1}$, which span the intersection hyperline. In this way, the continuity of the PwCA is ensured. After this final rotation the transformed $\vec{y}$ vector is the normal vector of each hyperplane that constitute the model (see Figure~\ref{fig:doubleConvex}\,d). The rotation angles of these final rotations are $r_{3,i}^-$ and $r_{3,i}^+$ for the hyperplanes on either side of the interface.

The last step, illustrated in Figure~\ref{fig:doubleConvex}\,e, is to compute the plane coefficients of the hyperplanes that constitute the model and of the interface that separates the two convex regions. They can be directly derived from each hyperplane's normal vector $\vec{\nu}$ (i.e.\ $\vec{y}$ for the model hyperplanes and $\vec{x}_1$ for the interface) and its origin $\vec{o}$.
\begin{equation}
    \vec{a} = \frac{1}{\| \vec{\nu}\|} \left[ -{\vec{o}}^\mathsf{T} \vec{\nu}, \nu_1, \ldots, \nu_n \right]^\mathsf{T}
    \label{eq:planeCoef}
\end{equation}

This series of transformation that defines the PwCA model is summarized in Algorithm~\ref{alg:transform}.

\begin{algorithm}
\SetAlgoLined
\SetKwInOut{Input}{input}\SetKwInOut{Output}{output}

\Input{parameters $\vec{r}_1, s_1, \vec{r}_{2,i}, s_{2,i}, r_{3,i}^-, r_{3,i}^+$}
\Output{hyperplane parameters $\vec{a}_i^-$, $\vec{a}_i^+$, interface $\vec{a}_\mathrm{ifc}$}

\SetKwFunction{rotate}{rotate}
\SetKwFunction{selectVector}{selectVector}
\SetKwFunction{getPlaneCoef}{getPlaneCoefficients}

\BlankLine
\tcp{The basis $\vec{B}$ is the aggregation of the basis vectors $[\vec{x}_1,\ldots,\vec{x}_{n-1},\vec{y}]$}
$\vec{B} \leftarrow$ identity matrix\tcp*[l]{orthonormal basis}
$\vec{o} \leftarrow$ $[0,\ldots,0]^\mathsf{T}$\tcp*[l]{origin of basis}

\BlankLine\tcp{Transformation 1}
$\vec{B} \leftarrow$ \rotate($\vec{B}$, $\vec{r}_1$)\;
$\vec{x}_1 \leftarrow$ \selectVector($\vec{B}$, 1)\;
$\vec{o} \leftarrow$ $\vec{o} + \vec{x}_1 \cdot s_1$\;
$\vec{a}_\mathrm{ifc} \leftarrow$ \getPlaneCoef($\vec{x}_1$, $\vec{o}$)\;

\BlankLine
\For{$i\leftarrow 1$ \KwTo $N_\mathrm{hyp}/2$}{
    \tcp{Transformation 2}
    $\vec{B}_i$ $\leftarrow$ \rotate($\vec{B}$, $\vec{r}_{2,i}$)\;
    $\vec{y} \leftarrow$ \selectVector($\vec{B}_i$, $n$)\;
    $\vec{o}_i$ $\leftarrow$ $\vec{o} + \vec{y} \cdot s_{2,i}$\;
    \BlankLine
    \tcp{Transformation 3, $-$ side of interface}
    $\vec{B}^\ast \leftarrow$ \rotate($\vec{B}_i$, $r_{3,i}^-$)\;
    $\vec{y} \leftarrow$ \selectVector($\vec{B}^\ast$, $n$)\;
    $\vec{a}_i^- \leftarrow$ \getPlaneCoef($\vec{y}$, $\vec{o}_i$)\;
    \BlankLine
    \tcp{Transformation 3, $+$ side of interface}
    $\vec{B}^\ast \leftarrow$ \rotate($\vec{B}_i$, $r_{3,i}^+$)\;
    $\vec{y} \leftarrow$ \selectVector($\vec{B}^\ast$, $n$)\;
    $\vec{a}_i^+ \leftarrow$ \getPlaneCoef($\vec{y}$, $\vec{o}_i$)\;
}

\caption{Convert rotation parameters to hyperplane parameters}
\label{alg:transform}
\end{algorithm}

In the end, $\vec{a}_i^-$ contains the parameters of the hyperplanes below the interface, i.e. in negative direction of the normal vector of the interface, $\vec{a}_i^+$ the parameters of hyperplanes above the interface and $\vec{a}_\mathrm{ifc}$ the parameters of the interface itself. Then the model estimate $\hat{y}$ is given by
\begin{equation}
    \hat{y} = \begin{cases}
        \max \{ y_i^-(\vec{x}) : i \in \{1,\ldots,N_\mathrm{hyp}/2\} \} \quad &\text{if} \quad 
        [1, \vec{x}, \hat{y}] \, \vec{a}_\mathrm{ifc}  \le 0 \\  
        
         \max \{ y_i^+(\vec{x}) : i \in \{1,\ldots,N_\mathrm{hyp}/2\} \} 
         \quad &\text{if} \quad 
         [1, \vec{x}, \hat{y}] \,  \vec{a}_\mathrm{ifc}  > 0 
    \end{cases}
\end{equation}
where $y_i^-$ are the $y$-values of the hyperplanes below the interface given by Eq.~\eqref{eq:yi} with the hyperplane parameters $\vec{a}_i^-$ and $y_i^+$ the equivalent above the interface.
        

The reformulated fitting problem is
\begin{equation}
    \argmin_{\vec{r}_1, s_1, \vec{r}_{2,i}, s_{2,i}, r_{3,i}^-, r_{3,i}^+} \sum_{m=1}^{N_\mathrm{data}} \big(\ \hat{y}(x'_{1,m},\ldots,x'_{n-1,m}) - y'_m \big)^2   .
    \label{eq:optDC}
\end{equation}

This is again an unconstrained non-linear optimization problem, which converges reliably given fairly good starting values. Just as the fitting algorithm for simple convex approximation, penalizing the distance between $y_i$ at the corners of the domain and $\max \{y'_m\}$ helps ensure that no hyperplane is dominated by the others on the whole domain.

\subsection{Starting values for piecewise-convex approximation}

To make the PwCA converge reliably, it is essential to compute good starting values for the hyperplanes that constitute the model. The key to compute good starting values is to observe that the intersection of the hyperplanes at the interface will be a convex function in $(n-1)$ dimension.

To initialize the model we start with an initial guess for the interface that separates the two regions $\vec{r}'_1, s'_1$.  By projecting points from the data set in the vicinity of the interface onto the interface and computing a simple convex approximation with the algorithm outlined in Section~\ref{sec:fitSC}, a good approximation of the intersection hyperlines $\vec{r}'_{2,i}, s'_{2,i}$ can be obtained. Then only one rotation parameter has to be estimated for each hyperplane $r_{3,i}^{ -\,\prime}, r_{3,i}^{+\,\prime }$ on either side of the interface to make the set of starting values complete. This can be done by using the optimization problem in Eq.~\eqref{eq:optDC} and fixing all parameters but the $r_{3,i}^-, r_{3,i}^+$.


\subsection{Translating a piecewise-convex approximation to MILP}


Assume that the PwCA is given by set of plane coefficients as described above: $\vec{a}_i^-$ contains the coefficients of the hyperplanes below the interface, i.e. in negative direction of the normal vector, $\vec{a}_i^+$ the coefficients above the interface and $\vec{a}_\mathrm{ifc}$ the coefficients of the interface itself. Then, the set of constraints that describes PwCA is

\begin{align}
     a_{\mathrm{ifc},0} + \sum_{j=1}^{n-1} a_{\mathrm{ifc},j} x_j + a_{\mathrm{ifc},n} y &\le M_\mathrm{t}^{+} t &
     \label{eq:ifc0}\\
     a_{\mathrm{ifc},0} + \sum_{j=1}^{n-1} a_{\mathrm{ifc},j} x_j + a_{\mathrm{ifc},n} y &\ge M_\mathrm{t}^{-} (1-t) & 
     \label{eq:ifc1}\\
     a_{i,0}^- + \sum_{j=1}^{n-1} a_{i,j}^- x_j + a_{i,n}^- y &\ge M_i^- t  &: \quad \forall i \in \{1,\ldots,N_\mathrm{hyp}/2\} 
     \label{eq:hyp0}\\
     a_{i,0}^+ + \sum_{j=1}^{n-1} a_{i,j}^+ x_j + a_{i,n}^+ y &\ge M_i^+ (1-t)   &: \quad \forall i \in \{1,\ldots,N_\mathrm{hyp}/2\}  
     \label{eq:hyp1}
\end{align}
where $x_1$,\ldots,$x_{n-1}$ and $y$ are the continuous variables of the optimization problem and $t$ is the binary toggle variable that is 0 if the point is below the interface and 1 if it is above the interface. Eq.~\eqref{eq:ifc0} and \eqref{eq:ifc1} set the toggle variable. Eq.~\eqref{eq:hyp0} and \eqref{eq:hyp1} represent the convex combination of the linear constraints, which constitute the PwCA and which are toggled according to $t$. (Note: For Eq.~\eqref{eq:hyp0} and \eqref{eq:hyp1} to work, the normal vector of the hyperplanes must point in positive $y$ direction, i.e. the $a_{i,n}$ component must have a positive sign.)

Translating a PwCA to MILP requires one binary variable and $N_\mathrm{hyp}+2$ big-M constraints. 
The big-M values are the maximum distances of points in the domain to the hyperplanes
\begin{align}
    M_t^+ &= \max_{\vec{p}} \vec{p}^\mathsf{T} \vec{a}_\mathrm{ifc} \\
    M_t^- &= \min_{\vec{p}} \vec{p}^\mathsf{T} \vec{a}_\mathrm{ifc} \\
    M_i^- &= \min_{\vec{p}} \vec{p}^\mathsf{T} \vec{a}_i^-  \quad :\quad \vec{p}^\mathsf{T} \vec{a}_\mathrm{ifc} \ge 0 \\
    M_i^+ &= \min_{\vec{p}} \vec{p}^\mathsf{T} \vec{a}_i^+  \quad :\quad \vec{p}^\mathsf{T} \vec{a}_\mathrm{ifc} \le 0
\end{align}
where
\begin{equation}
    \vec{p} = [1,x_1,\ldots, x_{n-1},y]^\mathsf{T} .
\end{equation}

Using this approach to translate concave functions to MILP takes \mbox{$N_\mathrm{hyp}/2+1$} binary variables and \mbox{$N_\mathrm{hyp}+2$} big-M constraints.

\subsection{Piecewise-convex approximation with three and more regions}

PwCA could be extended to more than two regions by using more than one interface. This would potentially make the approximation more accurate. Though, the  intersection hyperlines on two adjacent interfaces cannot be positioned independently, because they have to be connected by a hyperplane, which locks one degree of freedom. Consequently, the accuracy gained by adding additional interfaces will be diminished.
Furthermore, the main advantage of PwCA ---the highly efficient MILP formulation--- would suffer, because auxiliary binary variables for the additional interfaces would need to be introduced. 

PwCA is designed to fill the gap between simple convex approximation and approximation with simplices. To do this, one interface is sufficient. Considering the trade-off between accuracy, performance and complexity, we decided to not generalize PwCA in this way. 


\section{Evaluation}
\label{sec:results}
 
Piecewise convex approximation (PwCA) fills a gap between simple convex approximation and approximation with simplices. In this section we will first show that PwCA yields a significantly better approximation than simple convex approximation for functions that are curved in more than one direction. Then, we will show that PwCA yields significantly better performance than approximation with simplices in large optimization problems, because it introduces fewer binary variables.

As a benchmark, we use the multiplication of two variables $y = x_1 \cdot x_2$ in the range $x_1, x_2 \in [0, 1]$. This function has a positive curvature in the direction $x_1 + x_2$ and negative curvature orthogonal to that. Because of this curvature, simple convex approximation can be expected to perform poorly. Nevertheless, this function is a good example to showcase PwCA. For one, multiplication is, in a way, the most trivial kind of non-linearity and we need an efficient way to linearize it to solve non-linear problems with MILP. Additionally, there are various functions with similar shape that arise in engineering optimization problems.  

We use Matlab 2020b on a laptop with an Intel Core i5-10310U with 4 cores at 1.70\,GHz and 16\,GB DDR4 RAM. The optimization problems were solved with Matlab's \texttt{intlinprog} and the YALMIP toolbox \cite{lofberg2004}.

\subsection{Accuracy}

Here, we compare the accuracy of PwCA with simple convex approximation depending on the number of planes used for the approximation. The data set was generated by sampling the multiplication function on an equally spaced $100 \times 100$ grid. 

\begin{figure}
    \centering
    \makebox[\textwidth][c]{
            \includegraphics{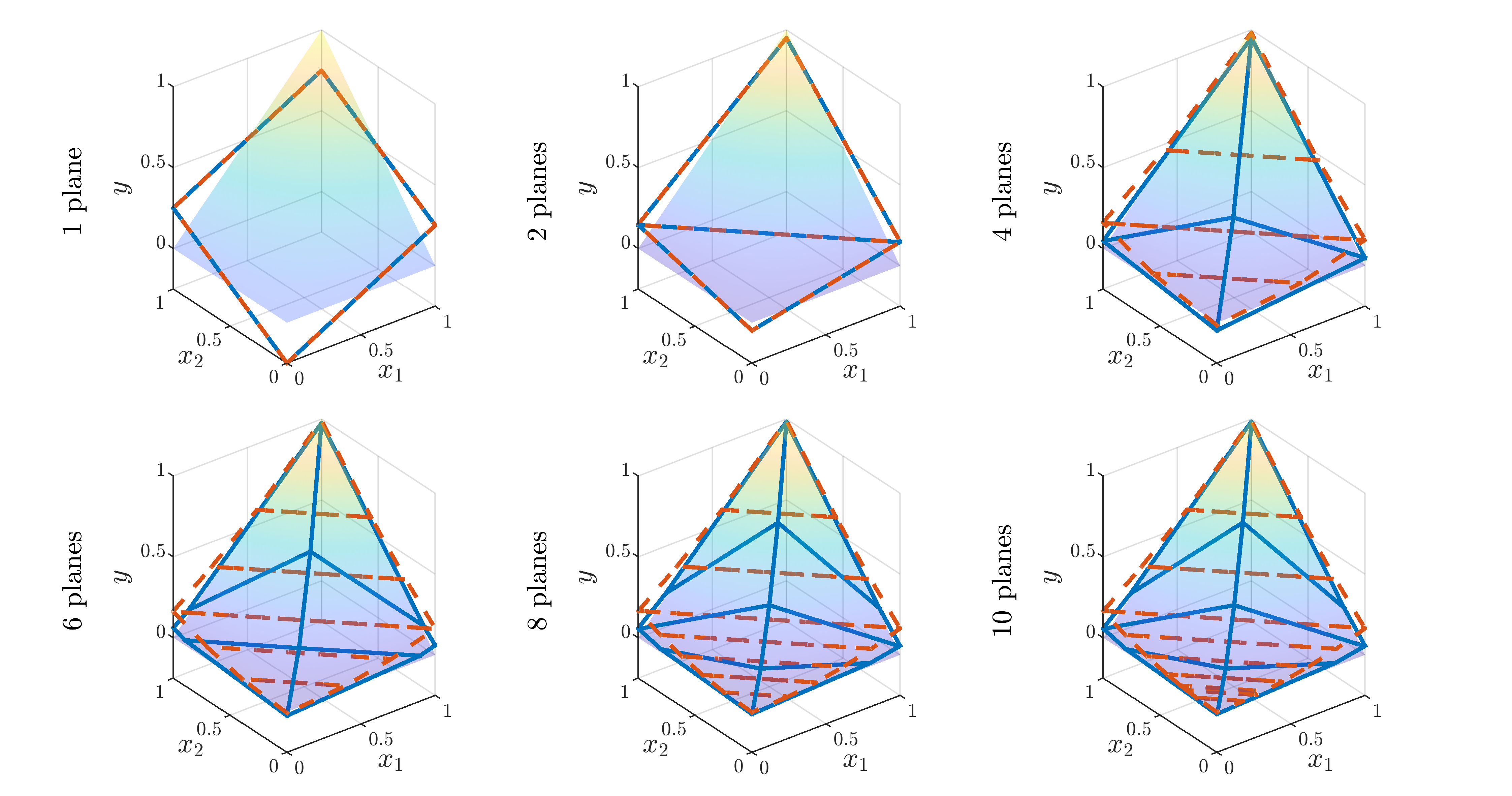}
    }
        \caption{Approximation of the function $y = x_1 \cdot x_2$ with different number of planes. Piecewise convex (blue solid), simple convex (red dashed).}
        \label{fig:comparison_models}
\end{figure}

\begin{figure}
    \centering
    \makebox[\textwidth][c]{
            \includegraphics{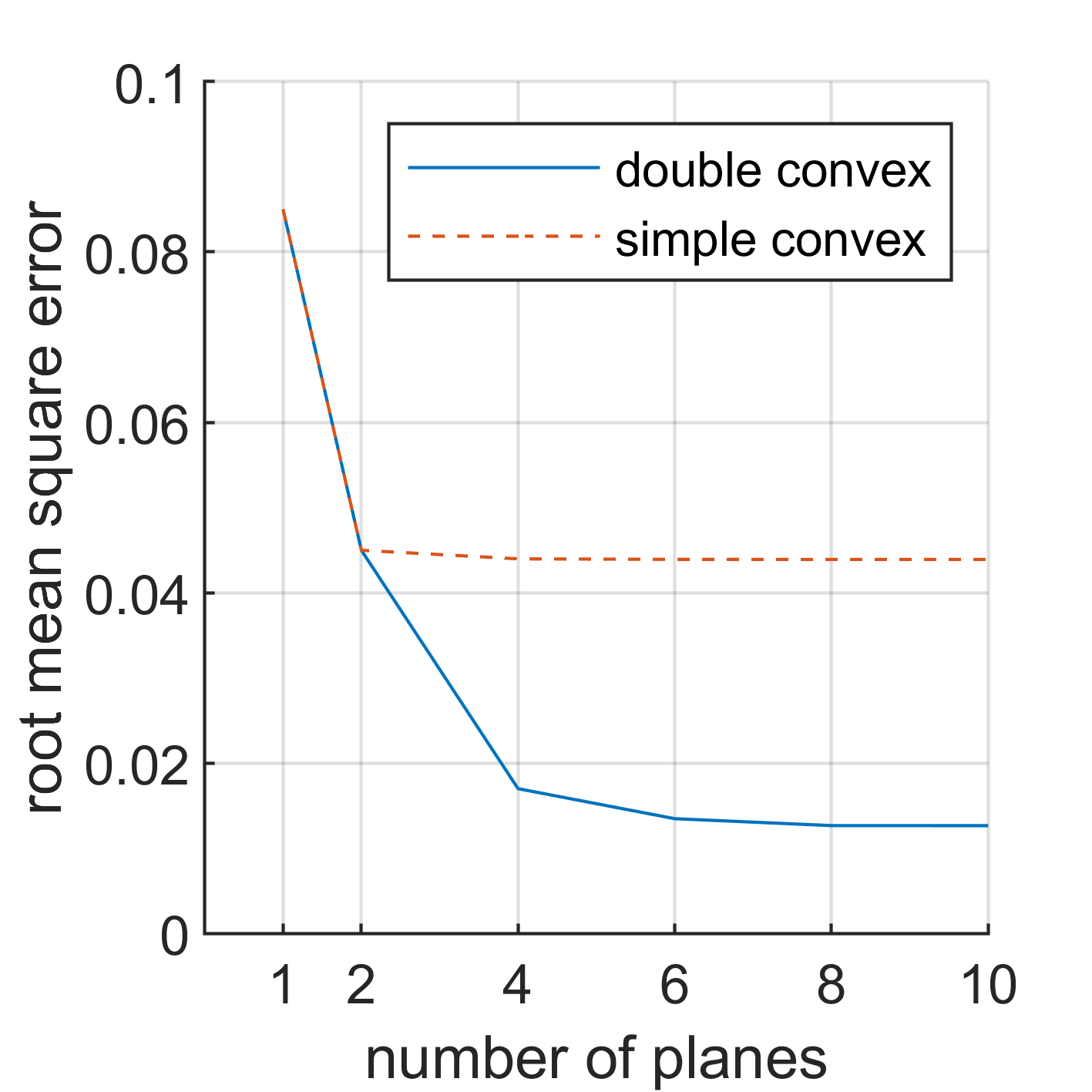}
    }
        \caption{Error of the approximations in Figure~\ref{fig:comparison_models}.}
        \label{fig:comparison_modelsRMSE}
\end{figure}

Figure~\ref{fig:comparison_models} shows the approximation of the multiplication function with 1 to 10~planes. The corresponding approximation error, measured by the root mean square error (RMSE), is shown in Figure~\ref{fig:comparison_modelsRMSE}. With the simple convex approximation anything more than two planes does not yield substantial improvements because then, the main curvature along the $x_1 + x_2$ direction is already well approximated and because simple convex approximation can ---\,by definition\,--- not model the other curvature.

PwCA, on the other hand, is more flexible. Starting at four planes, the PwCA splits the domain in half along the $x_1 + x_2$ direction and approximates both sides with a convex combination of two planes, which reduces the RMSE from 0.044 to 0.017. Adding two more planes, yields a slight improvement still. Then, also PwCA has reached its accuracy limit. If a more accurate approximation than this is required, one would  need to resort to simplex approximation.

Even though the accuracy of PwCA is limited, it is able to capture key characteristics of many non-linear functions. Depending on the application, this may be sufficient. Then, the main advantage of PwCA  ---\,the fast MILP solving time\,--- will come into effect, as we will see in the next subsection.

\subsection{MILP performance}

\begin{figure}
    \centering
    \includegraphics{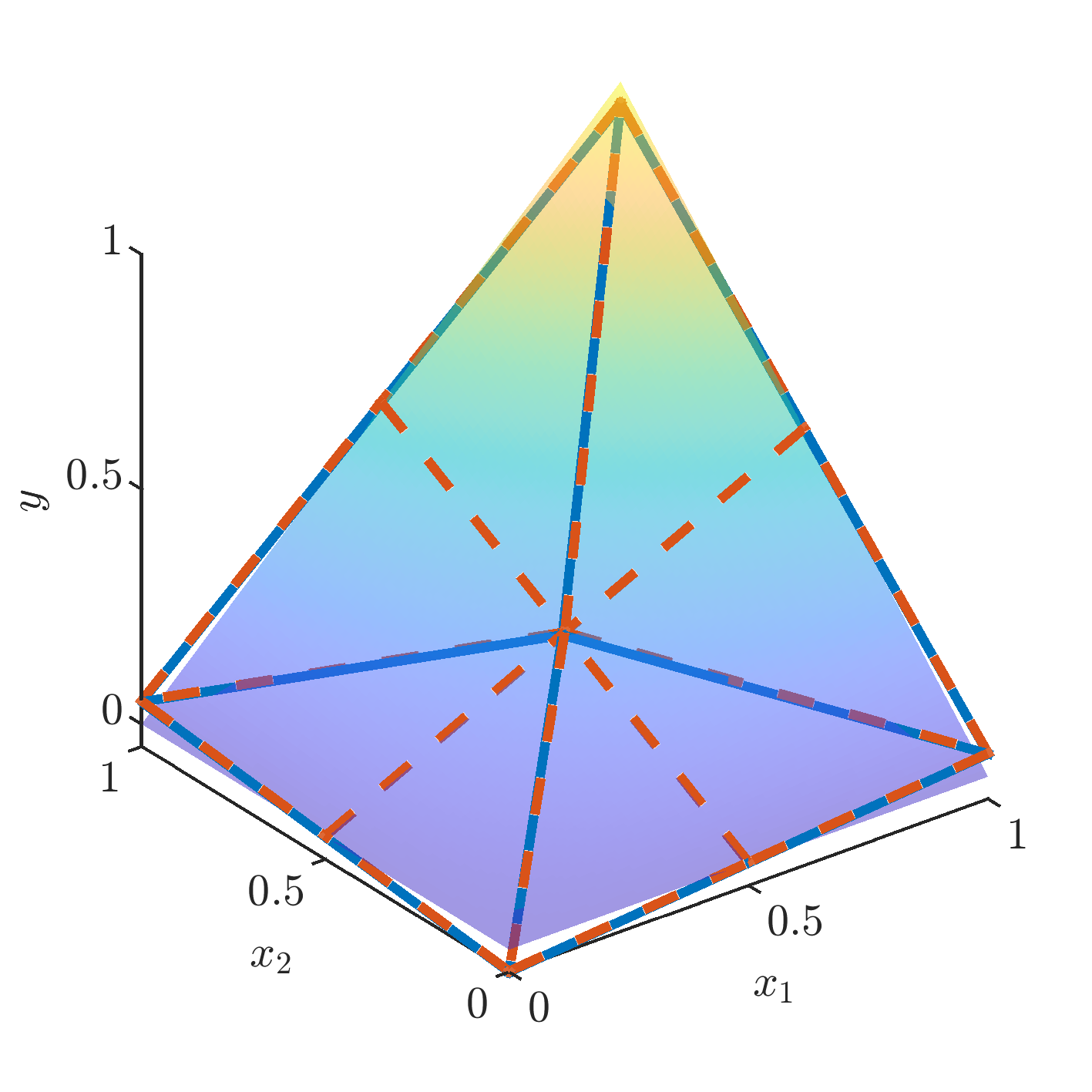}
    \caption{Illustration of the piecewise-convex approximation (blue solid) and the approximation using simplices (red dashed) for the benchmark.}
    \label{fig:models2D}
\end{figure}

To assess the performance of PwCA compared to approximation with simplices, the multiplication function was approximated using 4~planes. Then, an equivalent simplex approximation was set up, which required 8~simplices. Both these approximations achieved a root mean square error (RMSE) of 0.1304. The two approximations are shown in Figure~\ref{fig:models2D}. 


The translation of the PwCA to MILP was done with \mbox{Eq.~\eqref{eq:ifc0}--\eqref{eq:hyp1}} that were introduced in the previous section. For the translation of the simplex approximation, we chose three of the formulations introduced by Vielma et al. \cite{Vielma_2010}: convex combination (CC), multiple choice (MC) and the logarithmic convex combination (Log). Table~\ref{tab:my_label} shows the number of constraints and binary variables that are required for each model type and each MILP formulation. For a more detailed discussion on the number of variables and constraints in the simplex models, we refer to Table~1 in \cite{Vielma_2010}.

\begin{table}
    \centering
    \caption{Auxiliary variables and constraints for each approximation in Figure~\ref{fig:models2D}.}
    \begin{tabular}{l c c c}\hline
        approximation & binary & continuous & constraints \\ \hline
        piecewise convex & 1 & 0 & 5 \\
        simplex MC & 8 & 16 & 28 \\
        simplex CC & 8 & 9 & 13 \\
        simplex Log & 3 & 9 & 10 \\
        \hline
    \end{tabular}
    \label{tab:my_label}
\end{table}

For the benchmark, each approximation was used in a MILP problem, with a different number of query points $N$. In each run of the benchmark, $N$ random $[x'_{1,m}, x'_{2,m}]$ pairs were generated and the optimization problem was set up as 
\begin{align}
    \min &\sum_m y_m \\
    \text{s.t. } &y_m \ge \mathrm{PwLA}(x_{1,m} \cdot x_{2,m}) \\
                &x_{1,m} = x'_{1,m} \\
                &x_{2,m} = x'_{2,m} \\
                &y_m, x_{1,m}, x_{2,m} \in [0,1]
\end{align}
where PwLA is the MILP representation of the piecewise linear approximation of the product of $x_1$ and $x_2$.
For each query point, the MILP representation of the approximation must be replicated. Thus, the number of variables and constraints in the optimization problem increases linearly with $N$. 
The need to evaluate a model at various points arises, to name just one example, in unit commitment problems. If the model represents the operating behavior of a device, the model has to be replicated for each time step.  

\begin{figure}
    \centering
    \includegraphics{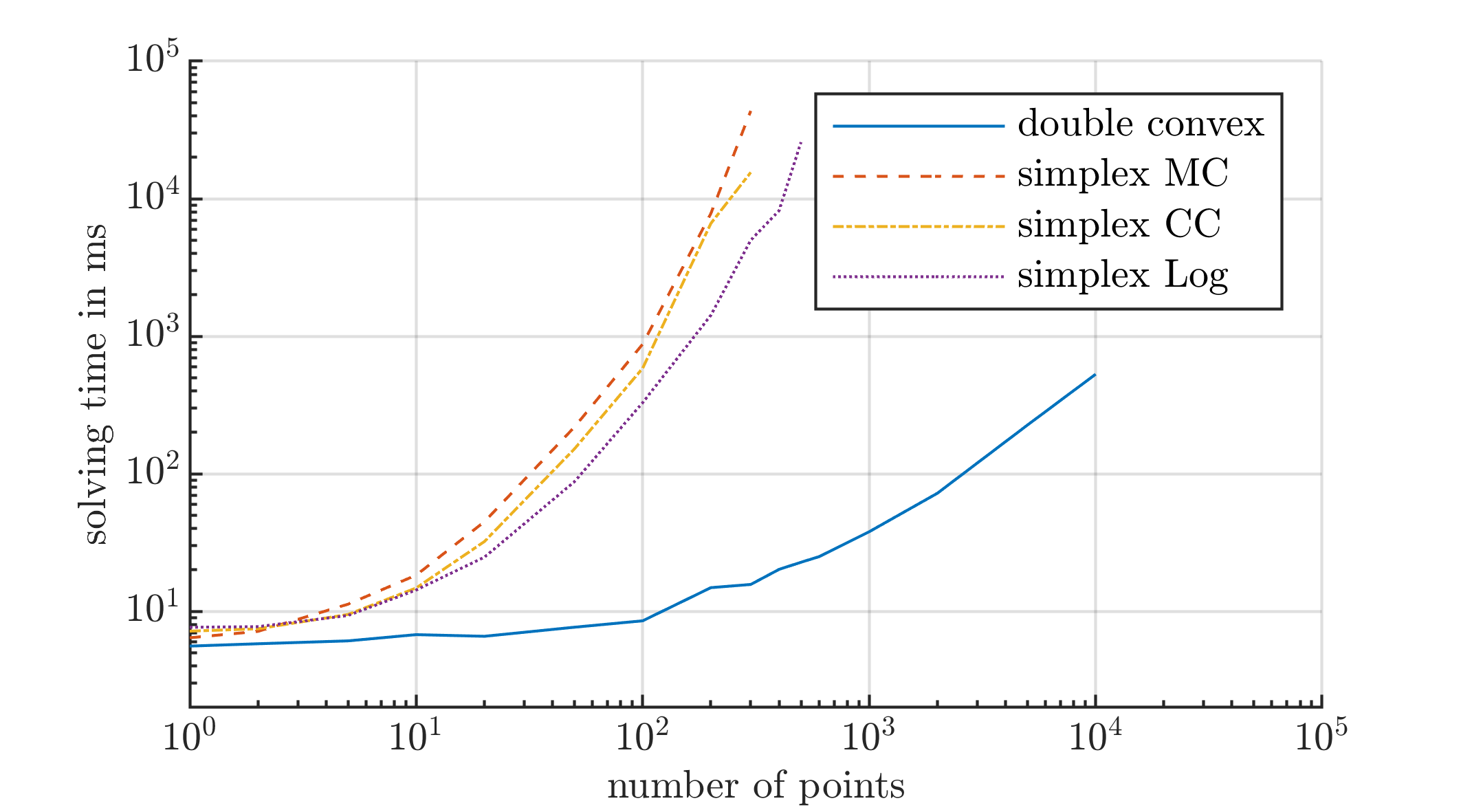}
    \caption{MILP performance of the PwCA and the simplex approximation with various MILP formulations in the test case.}
    \label{fig:performance2D}
\end{figure}

The results of the benchmark are illustrated in Figure~\ref{fig:performance2D}. We measured the performance in MILP problems with 1 to 10\,000 replications of the MILP representation of the approximations. To get a more accurate estimate of the solving time and eliminate random effects, each MILP problem was solved 10 times and the median solving time was computed. 

PwCA outperforms the simplex approximation by a large margin, regardless of the MILP formulation that was used to translate the simplex model to MILP. At 10~query points, PwCA performed 2~times faster than simplex approximation with Log-formulation (6.7\,ms vs.\ 14.3\,ms). At 300~query points, it even performed 39~times faster (8.5\,ms vs.\ 329\,ms), which is a consequence of the exponentially increasing complexity of MILP problems with the number of binary variables.

These results highlight the main advantage of PwCA: The MILP representation requires only one auxiliary binary variable. If the MILP representation is replicated multiple times in an optimization problem, e.g. at various time-steps, then the performance gain is significant.






\section{Conclusion}
\label{sec:conclusion}
We presented an algorithm to automatically compute a piecewise-linear approximation of a non-linear functions for MILP. The domain is split into two regions and each region is approximated by a convex combination of linear constraints. This piecewise-convex approximation (PwCA) fills a gap between simple convex approximation and approximation using simplices. 
PwCA can capture features of the non-linear function that simple convex approximation can not and it yields faster solving time than simplex approximation because only one auxiliary binary variable is required.

PwCA works especially well with functions that have both a direction with positive and with negative curvature. Then they provide a very accurate (in terms of deviation) approximation of the non-linear function, and at the same time they also perform very well in MILP. In this paper, we showcased this by considering the multiplication of two variables as an example. PwCA was both more accurate than simple convex approximation and yielded faster solving times than simplex approximation. Compared to simple convex approximation, the root mean square error was about 61\,\% lower (0.044 with simple convex vs. 0.017 with piecewise convex). 

Compared to approximation with simplices at the same accuracy, PwCA resulted in a 2~times faster solving time in a MILP problem with 10~instances and 39~times faster with 300~instances. This is a consequence of the exponentially increasing complexity of MILP problems with the number of binary variables. Since PwCA only requires one auxiliary binary variable (compared to 3 for an approximation with 8~simplices and Log-formulation) it scales considerably better for large MILP problems. 
This speed-up is especially useful for engineering optimization problems such as unit commitment, where the function that represents the operating behavior of a device has to be replicated at every time step. Then, the performance gain due to the reduced number of binary variables will be significant.

Research is under way to realize adaptive operation planning, where PwCA is used to approximate the operating behavior of devices in an energy system. PwCA fills a niche between simple convex approximation and simplex approximation, because it more accurate than simple convex approximation and yields faster solving times than simplex approximation. For this reason, PwCA will also be a useful tool for other applications where non-linear functions have to be approximated for MILP. 

\nomAcro[PwCA]{PwCA}{piecewise-convex approximation}
\nomAcro[PwLA]{PwLA}{piecewise-linear approximation}
\nomAcro[milp]{MILP}{mixed integer linear programming}
\nomAcro[milp]{MINLP}{mixed integer non-linear programming}
\nomAcro[rmse]{RMSE}{root mean square error}

\nomRoman[n]{$n$}{number of dimensions}
\nomRoman[n]{$N$}{number of points}
\nomRoman[nhyp]{$N_\mathrm{hyp}$}{number of hyperplanes}
\nomRoman[ndata]{$N_\mathrm{data}$}{number of points in data set}

\nomRoman[xj0]{$x_j$}{independent variable}
\nomRoman[x]{$\vec{x}$}{vector of independent variables}
\nomRoman[y]{$y$}{dependent variable}
\nomRoman[yest]{$\hat{y}$}{model estimate for dependent variable}

\nomRoman[xj]{$\vec{x}_j$}{vector in direction $x_j$}
\nomRoman[y]{$\vec{y}$}{vector in $y$-direction }

\nomRoman[a]{$\vec{a}$}{vector of hyperplane coefficients}
\nomRoman[a]{$a$}{hyperplane coefficient}

\nomRoman[r]{$r$}{rotation angle}
\nomRoman[s]{$s$}{shift distance}

\nomGreek[n]{$\vec{\nu}$}{normal vector of hyperplane}
\nomGreek[o]{$\vec{o}$}{origin of normal vector}

\nomSub[i]{$j$}{dimension index}
\nomSub[j]{$i$}{hyperplane index}
\nomSub[m]{$m$}{data set index}
\nomSub[ifc]{ifc}{interface}

\nomSuper[p]{$+$}{above the interface, i.e.\ in direction of the normal vector}
\nomSuper[m]{$-$}{below the interface, i.e.\ in negative direction of the normal vector}

\appendix


 \bibliographystyle{elsarticle-num} 
 \bibliography{references}


\printnomenclature




\end{document}